\documentclass[a4paper, 11pt]{amsart} 

\usepackage{verbatim} 
\usepackage{amssymb}
\usepackage{mathrsfs}
\usepackage[all]{xy}

\numberwithin{equation}{section}
\theoremstyle{plain}
\newtheorem{theorem}{Theorem}[section]
\newtheorem{proposition}[theorem]{Proposition}
\newtheorem{lemma}[theorem]{Lemma}
\newtheorem{corollary}[theorem]{Corollary}
\newtheorem{conjecture}[theorem]{Conjecture}

\theoremstyle{definition}
\newtheorem{definition}[theorem]{Definition}

\theoremstyle{remark}
\newtheorem{remark}[theorem]{Remark}

\newtheorem{example}[theorem]{Example}

\newcommand{\M}{\mathcal{M}}

\newcommand{\id}{\operatorname{id}}
\newcommand{\Ideal}{\operatorname{Ideal}}

\newcommand{\N}{\mathbb N}
\newcommand{\Z}{\mathbb Z}
\newcommand{\R}{\mathbb R}
\newcommand{\C}{\mathbb C}
\newcommand{\T}{\mathbb T}
\newcommand{\HilH}{\mathscr H}

\title{The ideal structure of reduced crossed products}
\author{Adam Sierakowski}
\address{Dept. of Mathematics \& Computer Science, University of Southern Denmark, Campusvej 55, DK-5230 Odense M, Denmark}
\email{adam@imada.sdu.dk}
\begin{document}

\begin{abstract}
Let $(A,G )$ be a C*-dynamical system with $G $ discrete. In this paper we investigate the ideal structure of the reduced crossed product C*-algebra and in particular we determine sufficient -- and in some cases also necessary -- conditions for $A$ to separate the ideals in $A\rtimes_r G $. When $A$ separates the ideals in $A\rtimes_r G $, then there is a one-to-one correspondence between the ideals in $A\rtimes_r G $ and the invariant ideals in $A$. We extend the concept of topological freeness and present a generalization of the Rokhlin property. Exactness properties of $(A, G )$ turns out to be crucial in these investigations.
\end{abstract}
\maketitle

In this paper we examine conditions on a C*-dynamical system $(A,G )$ with $ G $ discrete assuring that $A$ separates the ideals in the reduced crossed product $A\rtimes_r G $, i.e.\ when the (surjective) map $J\mapsto J\cap A$, from the ideals in $A\rtimes_r G $ into the invariant ideals in $A$, is injective.

Simplicity of $A\rtimes_r G $ obviously implies that $A$ separates the ideas in $A\rtimes_r G $. Some of the first results about simplicity of crossed products go back to works of Effros-Hahn \cite{EffHahn} and Zeller-Meier \cite{ZelMei}. For results restricted to abelian or Powers groups, we refer to \cite{OlePed3} and \cite{HarSkan}. Elliott showed in \cite{Ell} that $A\rtimes_r G $ is simple, provided that $A$ is an AF-algebra and the action is minimal and properly outer. Kishimoto showed in \cite{Kis:Auto} that the reduced crossed product of a C*-algebra $A$ by a discrete group $ G $ is simple if the action is minimal and fulfills the strong Connes spectrum condition.

Archbold and Spielberg generalized the result of Elliott by introducing topological freeness. The action of $ G $ on $A$ is called \emph{topologically free} if $\bigcap_{t\in F}\{x\in \widehat{A}: t.x\neq x\}$ is dense in the spectrum $\widehat{A}$ of $A$
(\footnote{
Here $\widehat{A}$ denotes the unitary equivalence classes of irreducible representations equipped with the not necessarily separated topology induced be the natural surjection onto the T$_0$ space $\mathrm{Prim}(A)$.
})
for any finite subset $F\subseteq G \setminus\{e\}$. They show in \cite{ArcSpi} that if the action is  minimal and topologically free, then the reduced crossed product $A\rtimes_r G $ is simple. For $(A, G )$ with $A$ abelian, the notion of topological freeness and proper outerness coincide, cf.\ \cite{ArcSpi}.

We say that an action of $ G $ on $A$ has the \emph{intersection property} if every non-trivial ideal in $A\rtimes_r G $ has a non-trivial intersection with $A$. The intersection property is a necessary condition to ensure that $A$ separates the ideals in $A\rtimes_r G $. {Kawamura} and Tomiyama showed in \cite{KawTom} that if $A$ is abelian and unital and $ G $ is amenable, then topological freeness of $ G $ acting on $A$ is equivalent to the intersection property. {There is also some recent work of Svensson and Tomiyama on a related ideal intersection property in the case where $A = C (X )$ and $G = \Z$, cf.\ \cite{Sven:Tom}. Archbold and Spielberg have made a generalization of the work by Kawamura and Tomiyama} to cover also non-amenable groups at the expense of replacing the intersection property with a stronger one, cf.\ Theorem \ref{lmt1}. Topological freeness is also connected to the property that $A$ is a maximal abelian sub-C*-algebra in $A\rtimes_r G $, cf.\ \cite[Proposition 4.14]{ZelMei}. Topological freeness alone is however neither sufficient nor necessary to ensure that $A$ separates the ideals in $A\rtimes_r G $. 

It was implicitly stated in the work by Renault, cf.\ \cite{Ren:groupoid}, that ``essential freeness" of $ G $ acting on $\widehat{A}$ might be enough to ensure that $A$ separates the ideals in $A\rtimes_r G $. The action of $ G $ on $\widehat{A}$ is called \emph{essentially free} provided that for every closed invariant subset $Y\subseteq \widehat{A}$, the subset of points in $Y$ with trivial isotropy is dense in $Y$. In the context of crossed products the claim of Renault can be phrased as follows:
\begin{conjecture}
\label{rencon}
{Let $(A, G)$ be a C*-dynamical system with $G$ discrete.} If the action of $ G $ on $\widehat{A}$ is essentially free then $A$ separates the ideals in $A\rtimes_r G $.
\end{conjecture}
{ Among other things} we show in this paper that, if the action of $G$ on $\widehat{A}$ is essentially free, then $A$ separates the ideals in $A\rtimes_r G $ if and only if the action of $G$ on $A$ is {'exact'. An action of $ G $ on $A$ is called \emph{exact} if every invariant ideal in $A$ induce a short exact sequence at the level of reduced crossed products ({see} Definition \ref{deo2}).} The necessity follows from Theorem \ref{thetheorem}, and the sufficiency follows from Theorem \ref{tht1}. It is not at all clear, and is an important point, if we have automatic exactness of the action under the assumption of essential freeness, compare Remark \ref{rem:Gromov}. 

For a separable C*-dynamical system $(A, G )$, i.e.\ with $ G $ countable and $A$ separable, Effros and Hahn showed that $A$ separates the ideals in $A\rtimes_r G $ when $A$ is abelian and when $ G $ is amenable and acts freely on $\widehat{A}$, cf.\ \cite{EffHahn}. 

Zeller-Meier made a generalization of that result to GCR algebras. He showed in \cite{ZelMei} that for a separable C*-dynamical $(A, G )$ with $ G $ discrete, $A$ separates the ideals in $A\rtimes_r G $ provided $A$ is GCR and the action of $ G $ on $\widehat{A}$ (which is for GCR algebras naturally isomorphic to $\mathrm{Prim}(A)$) is free and regular. The action is called \emph{regular} if every $ G $-ergodic measure on $\widehat{A}$ is supported in a single $ G $-orbit. We remark that the conditions, in the work by Effros-Hahn and by Zeller-Meier, ensuring that $A$ separates the ideals in $A\rtimes_r G $ also imply that the canonical maps $\pi^{A/I}: (A/I)\rtimes  G \to (A/I)\rtimes_r G $ are isomorphisms {for all closed, two-sided $G$-invariant ideals $I$ of $A$.} The latter property is equivalent to exactness of the action and injectivity of $\pi^A$ (by exactness of the full crossed product functor and application of the 5-lemma). We don't know if exactness of the action of $G$ on $A$ {follows from the injectivity of $\pi^A$.} A counterexample needs an improvement of Gromov's construction in \cite{Gromov:exact}.

Renault considered in \cite{Ren:fixed} the groupoid crossed product C*-algebras. His work contains some generalizations of the results in \cite{EffHahn, ZelMei}, where freeness of $ G $ acting on $\widehat{A}$ is replaced with essential freeness.

In the first part of this article we show that $A$ separates the ideals in $A\rtimes_r G $ provided that the action of a $ G $ on $A$ is exact and the action of $G$ on $\widehat{A}$ is essentially free.	The exactness of the action is a \emph{necessary condition} to ensure that $A$ separates the ideals in $A\rtimes_r G $. Certainly, the exactness of the action alone does not imply that $A$ separates the ideals of $A\rtimes _r G$ (even in case of abelian $G$).

In the second part of the article we present another way of ensuring that $A$ separates the ideals in $A\rtimes_r G $ by means of a Rokhlin type property. We use this result to make a generalized version of the Rokhlin property, which we call the residual Rokhlin* property. We show that $A$ separates the ideals in $A\rtimes_r G $ provided that the action of $ G $ on $A$ is exact and satisfy the residual Rokhlin* property. 

In the case of \emph{abelian} $A$ and countable $G$ we have 
the following properties are equivalent, c.f.\ Corollary \ref{cor.additional.remark} and Remark \ref{twotopfree}:
\begin{itemize}
\item[(i)] For every invariant ideal  $I$ in $A$, the action of $G$ on $A/I$ is topologically free in the sense of Archbold and Spielberg.
\item[(ii)] The action of $G$ on $\widehat{A}$ is essentially  free in the sense  of  Renault.
\item[(iii)] The action of $G$ on $A$ satisfy the residual Rokhlin* property.
\item[(iv)] For every invariant ideal $I\not=A$ and every $t\neq e$ the automorphism $[a]\mapsto [t.a]$ on $A/I$ is properly outer in the sense of Elliott.
\end{itemize}

My profound thanks go to Professor E.\ Kirchberg for his extensive feedback and help on this paper --- including the ideas behind the residual Rokhlin* property --- and Professor M.\ R{\o}rdam for his guidance. I am also very grateful for the suggestions, contributions, and bug reports offered by many people, including Professor G.\ A.\ Elliott and Professor N.\ C.\ Phillips.
\setcounter{section}{0}
\section{Condition expectation}
\subsection{Notation}
Let $e$ denote the unit of a group. Ideals are always assumed to be closed and two-sided. Let $\M(A)$ denote the multiplier algebra of $A$ and let $A\rtimes_r G $ (resp. $A\rtimes G $) denote the reduced (resp. the full) crossed product. The reduced norm will be written as $\|\cdot\|_\lambda$. Let $\mathcal{I}(A)$ denote the set of ideals in a C*-algebra $A$ and let $\textrm{Ideal}_A[S]$, or simply $\textrm{Ideal}[S]$, be the smallest ideal in $A$ generated by $S\subseteq A$. Fix now a C*-dynamical system $(A, G )$ with $ G $ discrete. We will suppress the canonical inclusion map $A\subseteq A\rtimes_r G $. We let $\mathcal{I}(A)^G$ denote the set of all \emph{invariant} (or more precisely all \emph{$G$-invariant}) ideals in $A$ and let $\textrm{Ideal}_A[S]^ G $, or simply $\textrm{Ideal}[S]^ G $, be the smallest invariant ideal in $A$ generated by $S\subseteq A$. Elements in $C_c( G ,A)$ will be written as finite sums $\sum_{t\in F}a_tu_t$, where $u_t$ is, in a canonical way, a unitary element in $\M(A\rtimes_r G )$.

We now mention some of the well know (or easy to prove) results concerning reduced crossed products. { Let $(A, G )$ be a C*-dynamical system with $ G $ discrete.} The map $E\colon C_c( G ,A)\to A, \sum_{t\in F}a_tu_t\mapsto a_e$ extends by continuity to a faithful positive conditional expectation $E:A\rtimes_r G \to A$. Sometimes one writes $E_A$ instead of $E$. For every $I\in \mathcal{I}(A)^ G $ the natural maps in the short exact sequence
\begin{displaymath}
\xymatrix{
0 \ar[r] & I \ar[r]^-{\iota} & A \ar[r]^{\rho} & A/I \ar[r] & 0}
\end{displaymath}
extend, in a canonical way, to maps at the level of reduced crossed products and give the following commutative diagram (cf.\ \cite[Remark 7.14]{Will:cross})	
\begin{eqnarray}
\label{diagram}
\xymatrix{0 \ar[r] & I\rtimes_r G  \ar[r]^-{\iota \rtimes_r \id} \ar[d]^-{E_I} & A\rtimes_r G  \ar[r]^{\rho \rtimes_r \id} \ar[d]^-{E_A}& A/I\rtimes_r G \ar[r] \ar[d]^-{E_{A/I}}& 0\\
0 \ar[r] & I \ar[r]^-{\iota} & A \ar[r]^{\rho} & A/I \ar[r] & 0}
\end{eqnarray}
Note that the sequence at the level of reduced crossed products is only short exact provided that $\ker(\rho\rtimes_r \id )\subseteq I\rtimes_r G $. This inclusion does not hold in general, cf.\ Remark \ref{rem:Gromov}. However when $ G $ is exact  we do have that $\ker(\rho\rtimes_r \id )\subseteq I\rtimes_r G $, cf.\ \cite{KirWas}. With $I\in \mathcal{I}(A)^ G $ we have the identities
$$I\rtimes_r G =\textrm{Ideal}_{A\rtimes_r G}[I], \quad  (I\rtimes_r G )\cap A=I=E_A(I\rtimes_r G )\,,$$
and for $J\in \mathcal{I}(A\rtimes_r G )$ it follows that $J\cap A\in \mathcal{I}(A)^ G $ is a subset of $E(J)$. For an introduction to crossed products we refer to \cite{BroOza, Will:cross}.
\subsection{Ideals in $A\rtimes_r G $}
Fix a C*-dynamical system $(A, G )$ with $ G $ discrete. The question concerning when the algebra separates the ideals in the reduced crossed product, i.e.\ when the map
$$\mathcal{I}(A\rtimes_r G )\to \mathcal{I}(A)^ G : J\mapsto J\cap A$$
in injective, can be reformulated in several ways. Note that this map is automatically surjective (using that $(I\rtimes_r G ) \cap A=I$ for $I\in \mathcal{I}(A)^ G $). Hence, if $A$ separates the ideals in $A\rtimes_r G $, then there is a ono-to-one correspondence between ideals in the crossed product $A\rtimes_r G $ and the invariant ideals in $A$. 
\begin{proposition}
\label{lmn2}
{Given a C*-dynamical system $(A, G )$ with $ G $ discrete, the following properties are equivalent:}
\begin{itemize}
\item[(i)] For every $J\in \mathcal{I}(A\rtimes_r G )$ we have that $J=\Ideal_{A\rtimes_r G}[J\cap A]$.
\item[(ii)] The map $\mathcal{I}(A\rtimes_r G )\to \mathcal{I}(A)^ G : J\mapsto J\cap A$ is injective.
\item[(iii)] The map $\mathcal{I}(A)^ G \to \mathcal{I}(A\rtimes_r G ): I\mapsto I\rtimes_r G $ is surjective.
\item[(iv)] For every $J\in \mathcal{I}(A\rtimes_r G )$ we have that $J=\Ideal_{A\rtimes_r G}[E(J)]$.
\end{itemize}
\end{proposition}
\begin{proof}
$(i)\Rightarrow (ii)$ {Fix} ideals $J_1,J_2\in \mathcal{I}(A\rtimes_r G )$ having the same intersection with $A$, one gets $$J_1=\textrm{Ideal}[J_1\cap A]=\textrm{Ideal}[J_2\cap A]=J_2.$$

$(ii)\Rightarrow (iii)$ {Let} $J\in \mathcal{I}(A\rtimes_r G )$. Using the general fact that $$((J\cap A)\rtimes_r  G )\cap A=J\cap A$$ together with $(ii)$ it follows that $(J\cap A)\rtimes_r  G =J$.

$(iii)\Rightarrow (iv)$ {Let} $J\in \mathcal{I}(A\rtimes_r G )$. From $(iii)$ it follows that $J=I\rtimes_r  G $ for some  $I\in \mathcal{I}(A)^ G $. We get that 
$$J=I\rtimes_r  G =\textrm{Ideal}[I]=\textrm{Ideal}[E_A(I\rtimes_r  G )]=\textrm{Ideal}[E_A(J)].$$

$(iv)\Rightarrow (i)$ {Let}
$J\in \mathcal{I}(A\rtimes_r G )$. By $(iv)$ we have that $$J\cap A\subseteq E_A(J)\subseteq \textrm{Ideal}[E_A(J)]\cap A=J\cap A.$$ Using $(iv)$ once more we obtain that $J=\textrm{Ideal}[E_A(J)]=\textrm{Ideal}[J\cap A]$.
\end{proof}
The equality in part $(i)$ of Proposition \ref{lmn2} has one trivial inclusion. This is not the case when we consider the equality in $(iv)$. Here we have two non-trivial inclusion. However, as we will see, one of the inclusions corresponds to exactness of the action (cf.\ Definition \ref{deo2}), and hence it is automatic provided that $ G $ is exact.
\subsection{Exactness and the residual intersection property}
Let us now consider the two inclusions corresponding to the last equality in Proposition \ref{lmn2} separately.
\begin{definition}
\label{deo2}
Let $(A, G )$ be a C*-dynamical system. We say that the action (of $ G $ on $A$) is \emph{exact} if every invariant ideal $I$ in $A$ induces a short exact sequence 
\begin{displaymath}
\xymatrix{0 \ar[r] & I\rtimes_r G  \ar[r]^-{\iota \rtimes_r \id} & A\rtimes_r G  \ar[r]^{\rho \rtimes_r \id} & A/I\rtimes_r G \ar[r] & 0}
\end{displaymath}
at the level of reduced crossed products. A group $G$ is called \emph{exact} if any action of $G$ on any C*-algebra is exact, cf.\ \cite{KirWas}. 
\end{definition}
\begin{proposition}
\label{tho1}
{ Let $(A, G )$ be a C*-dynamical system with $ G $ discrete.} Then the following properties are equivalent
\begin{itemize}
\item[(i)] For every $J\in \mathcal{I}(A\rtimes_r G )$ we have that $J\subseteq \Ideal_{A\rtimes_r G}[E(J)]$.
\item[(ii)] For every $x\in (A\rtimes_r G )^+$ we have that $x\in \Ideal_{A\rtimes_r G}[E(x)]$.
\item[(iii)] The action of $ G $ on $A$ is exact.
\end{itemize}
\end{proposition}
\begin{proof}
$(ii)\Rightarrow (i)$ From 
$$J^+=\bigcup_{x\in J^+}\{x\}\subseteq \bigcup_{x\in J^+}\textrm{Ideal}[E(x)]\ \subseteq \ \textrm{Ideal}[E(J)]$$ 
it follows that $J=\textrm{Ideal}[J^+]\subseteq \textrm{Ideal}[E(J)]$.

$(i)\Rightarrow (iii)$ {Let} $I\in \mathcal{I}(A)^ G $. It is sufficient to verify the inclusion $\ker (\rho \rtimes_r \id )\subseteq I\rtimes_r G $ for the commutative diagram (\ref{diagram})
\begin{displaymath}
\xymatrix{0 \ar[r] & I\rtimes_r G  \ar[r]^-{\iota \rtimes_r \id} \ar[d]^-{E_I} & A\rtimes_r G  \ar[r]^{\rho \rtimes_r \id} \ar[d]^-{E_A}& A/I\rtimes_r G \ar[r] \ar[d]^-{E_{A/I}}& 0\\
0 \ar[r] & I \ar[r]^-{\iota} & A \ar[r]^{\rho} & A/I \ar[r] & 0.}
\end{displaymath}
With $J:=\ker(\rho \rtimes_r \id )$ we get that $\rho(E_A(J))=E_{A/I}(\rho \rtimes_r \id (J))=0$. This implies that $E_A(J)\subseteq \ker \rho=I$ and hence that $$J\subseteq \textrm{Ideal}[E_A(J)]\subseteq \textrm{Ideal}[I]=I\rtimes_r G .$$

$(iii)\Rightarrow (ii)$ {Let} $x\in (A\rtimes_r G )^+$. Set $I:=\textrm{Ideal}_A[E_A(x)]^ G $ and consider the two, using $(iii)$, short exact sequences in the commutative diagram (\ref{diagram})
\begin{displaymath}
\xymatrix{0 \ar[r] & I\rtimes_r G  \ar[r]^-{\iota \rtimes_r \id} \ar[d]^-{E_I} & A\rtimes_r G  \ar[r]^{\rho \rtimes_r \id} \ar[d]^-{E_A}& A/I\rtimes_r G \ar[r] \ar[d]^-{E_{A/I}}& 0\\
0 \ar[r] & I \ar[r]^-{\iota} & A \ar[r]^{\rho} & A/I \ar[r] & 0.}
\end{displaymath}
As $E_A(x)\in I$, it follows that $E_{A/I}(\rho \rtimes_r \id (x))=\rho(E_A(x))=0$. As $E_{A/I}$ is faithful on positive elements, it follows that $\rho \rtimes_r \id (x)=0$ and hence
$$x\in \ker(\rho \rtimes_r \id )=I\rtimes_r G =\textrm{Ideal}[I]=\textrm{Ideal}[\textrm{Ideal}_A[E_A(x)]^ G ]= \textrm{Ideal}[E_A(x)].$$
\end{proof}
Let us present the first application of Proposition \ref{tho1}. In general it is unknown whenever a family of invariant ideals $(I_i)$ fulfills that $\bigcap_i(I_i\rtimes_r G )=(\bigcap_i I_i)\rtimes_r G $. However when the action is exact the equality is easily shown.
\begin{proposition}
\label{coo2}
Let $(A, G )$ be a C*-dynamical system with $ G $ discrete. If the group $ G $ (or the action) is exact, then every family of ideals $(I_i)$ in $\mathcal{I}(A)^ G $ fulfills the identity
$$\bigcap_i(I_i\rtimes_r G )=(\bigcap_i I_i)\rtimes_r G .$$
\end{proposition}
\begin{proof}
Note that $E(\bigcap_i(I_i\rtimes_r G ))\subseteq \bigcap_i I_i$. By Proposition \ref{tho1} we get the inclusion  $\bigcap_i(I_i\rtimes_r G )\subseteq \textrm{Ideal}[E(\bigcap_i(I_i\rtimes_r G ))]$ and hence
\begin{eqnarray*}
(\bigcap_i I_i)\rtimes_r G &\subseteq& \bigcap_i(I_i\rtimes_r G )\ \subseteq \ \textrm{Ideal}[E(\bigcap_i(I_i\rtimes_r G ))]\\
&\subseteq& \textrm{Ideal}[\bigcap_i I_i]= (\bigcap_i I_i)\rtimes_r G .
\end{eqnarray*}
\end{proof}

\begin{remark}
A classical question of non-commutative harmonic analysis is how the irreducible representations
of $A\rtimes G$ or at least of $A\rtimes _r G$ look like. {Since any proper, closed, two-sided ideal of a C*-algebra is the intersection of the primitive ideals containing it, we can conclude from Proposition 1.1 
and the semi-continuity} property in Proposition \ref{coo2} that at least the following holds:
{\it If the action of $G$ on $A$ is exact, then every irreducible representation $D$ of $A\rtimes_r G$ has kernel $J_D=(A\cap J_D)\rtimes_r G$ if and only if $A$ separates the ideals of $A\rtimes_r G$.}
\end{remark}	
The second inclusion in part $(iv)$ of Proposition \ref{lmn2} is closely related to the way ideals in the reduced crossed product intersect the original algebra.
\begin{definition}
\label{deo3}
Let $(A, G )$ be a C*-dynamical system with $ G $ discrete. We say that the action (of $ G $ on $A$) has the \emph{intersection property} if every non-zero ideal in $A\rtimes_r G $ has a non-zero intersection with $A$. If the induced action of $ G $ on $A/I$ has the intersection property for every invariant ideal $I$ in $A$, we say that action (of $ G $ on $A$) has the \emph{residual intersection property}.
\end{definition}
\begin{remark}
\label{lmp0}
Let $(A, G )$ be a C*-dynamical system with $ G $ discrete. The intersection property of the action is equivalent to the requirement that every representation of $A\rtimes_r G $, which is faithful {on $A$, is itself faithful.} Hence a way to ensure that the action has the intersection property is to show that every representation of $A\rtimes G $, which is faithful on $A$, weakly contains the regular representation of $(A, G )$. (For a proof of the first claim consider restricting the representations of $A\rtimes_r G $. The second part uses the fact that the surjection $\pi\colon A\rtimes  G \to  A\rtimes_r  G $ is the identity map on $C_c( G ,A)$.)
\end{remark}	
\begin{proposition}
\label{tho2}
{ Let $(A, G )$ be a C*-dynamical system with $ G $ discrete.} Then the following properties are equivalent
\begin{itemize}
\item[(i)] For every $J\in \mathcal{I}(A\rtimes_r G )$ we have that $J\supseteq \Ideal_{A\rtimes_r G}[E(J)]$.
\item[(ii)] For every $x\in (A\rtimes_r G )^+$ we have that $E(x)\in \Ideal_{A\rtimes_r G}[x]$.
\item[(iii)] The action {satisfies} the residual intersection property and for every $J\in \mathcal{I}(A\rtimes_r G )$ the intersection $(\rho\rtimes_r \id)(J)\cap A/I$ is zero, where $I:=J\cap A$ and $\rho\rtimes_r \id $ comes from (\ref{diagram}).
\end{itemize}
\end{proposition}
\begin{proof}
$(i)\Rightarrow (iii)$: First we show the residual intersection property. Take $I\in \mathcal{I}(A)^ G $ and $J\in \mathcal{I}(A/I\rtimes_r  G )$ and assume that $J\cap A/I=0$. We show that $E_{A/I}(J)=0$, which is equivalent to $J=0$. Set $J_1:=(\rho\rtimes_r \id )^{-1}(J)$ where $\rho\rtimes_r \id $ comes from (\ref{diagram}). As $J_1\in \mathcal{I}(A\rtimes_r G )$ it follows from $(i)$ that $E_A(J_1)\subseteq J_1\cap A$ and hence
$$E_{A/I}(J)=E_{A/I}((\rho\rtimes_r \id) (J_1))=\rho(E_A(J_1))\subseteq \rho(J_1\cap A)\subseteq J\cap A/I=0.$$

For the second part take $J\in \mathcal{I}(A\rtimes_r  G )$ and $x\in(\rho\rtimes_r \id )(J)\cap A/I$ where $I:=J\cap A$ and $\rho\rtimes_r \id $ comes from (\ref{diagram}). We show $x=0$. Find $j\in J$ such that $x=(\rho\rtimes_r \id) (j)\in A/I$. Using $(i)$ we have that $E_A(J)\subseteq J\cap A$ and hence
$$x=E_{A/I}(x)=E_{A/I}((\rho\rtimes_r \id) (j))=\rho(E_A(j))\in \rho(E_A(J))\subseteq\rho(J\cap A)=0.$$
	
$(iii)\Rightarrow (ii)$ {Let} $x\in (A\rtimes_r G )^+$. Set $J:=\textrm{Ideal}[x]$ and $I:=J\cap A$. Using $(iii)$ on $J\in \mathcal{I}(A\rtimes_r G )$ it follows that
$$(\rho\rtimes_r \id )(J)\cap A/I=0,$$
with the surjection $\rho\rtimes_r \id $ coming from (\ref{diagram}). The residual intersection property implies that the ideal $(\rho\rtimes_r \id)(J)=0$. Using the diagram (\ref{diagram}) we now have that $\rho(E_A(J))=E_{A/I}(\rho\rtimes_r \id (J))=0$ and hence
$$E_A(x)\in E_A(J)\subseteq \ker \rho=J\cap A\subseteq J=\textrm{Ideal}[x].$$

$(ii)\Rightarrow (i)$ {Let} $J\in \mathcal{I}(A\rtimes_r G )$. Then 
$$E_A(J^+)=\bigcup_{x\in J^+}\{E_A(x)\}\subseteq \bigcup_{x\in J^+}\textrm{Ideal}[x]\subseteq J.$$
As every element in $J$ is a linear combination of positive elements in $J$ and $E_A$ is linear, it follows that $E_A(J)\subseteq J$ and hence also $\textrm{Ideal}[E_A(J)]\subseteq J$.
\end{proof}
Using an additional observation, contained in Lemma \ref{lemo3} below, we obtain a new characterization of when $A$ separates the ideals in $A\rtimes_r G $.
\begin{lemma}
\label{lemo3}
{ Let $(A, G )$ be a C*-dynamical system with $ G $ discrete.} Suppose the action of $ G $ on $A$ is exact. Then for every $J\in \mathcal{I}(A\rtimes_r G )$ the intersection $(\rho\rtimes_r \id) (J)\cap A/I$ is zero, where $I:=J\cap A$ and $\rho\rtimes_r \id $ comes from (\ref{diagram}).
\end{lemma}
\begin{proof}
For a given $J\in \mathcal{I}(A\rtimes_r G )$ set $I:=J\cap A$. Using that the action is exact we have the short exact sequence
\begin{displaymath}
\xymatrix{0 \ar[r] & I\rtimes_r G  \ar[r]^-{\iota \rtimes_r \id} & A\rtimes_r G  \ar[r]^{\rho \rtimes_r \id} & A/I\rtimes_r G \ar[r] & 0}
\end{displaymath}
and can therefore identify $(\rho\rtimes_r \id) (A\rtimes_r G )$ with $A\rtimes_r G /L$, where $L:=I\rtimes_r G $. We obtain the identities
$$J\cap A=I=L\cap A,$$
$$(\rho\rtimes_r \id )(J)=J/L, \ \ \ A/I=A/(L\cap A)=(A+L)/L.$$
Assume $(\rho\rtimes_r \id) (J)\cap A/I\neq 0$. {Then there exists $j\in J$ and $a\in A$ such that
$$j+L=a+L\neq L.$$
As $L\subseteq J$ it follows that $a\in J$ and hence $a\in J\cap A=I\subseteq L$. But this implies $a+L=L$} and we get a contradiction. Hence $(\rho\rtimes_r \id) (J)\cap A/I=0$.
\end{proof}
\begin{theorem}
\label{thetheorem}
{ Let $(A, G )$ be a C*-dynamical system with $ G $ discrete.} Then the following properties are equivalent. 
\begin{itemize}
\item[(i)] $A$ separates the ideals in $A\rtimes_r  G $. 
\item[(ii)] The action is exact and for every $x\in (A\rtimes_r G )^+: E(x)\in \Ideal_{A\rtimes_r G}[x]$.
\item[(iii)] The action is exact and satisfies the residual intersection property.
\end{itemize}
\end{theorem}
\begin{proof}
Combine Proposition \ref{lmn2}, Proposition \ref{tho1} and Proposition \ref{tho2} with Lemma \ref{lemo3}.
\end{proof}
Using Proposition \ref{coo2} one can slightly improve the last part of Theorem \ref{thetheorem} in the following sense:
\begin{corollary}
\label{coo3}
{ Let $(A, G )$ be a C*-dynamical system with $ G $ discrete.} Then the following properties are equivalent.
\begin{itemize}
\item[(i)] $A$ separates the ideals in $A\rtimes_r  G $. 
\item[(ii)] The action is exact and the intersection $J\cap A/I$ is non-zero for every $I\in \mathcal{I}(A)^ G $ and for every non-zero primitive ideal $J\in \mathcal{I}(A/I\rtimes_r G )$. 
\end{itemize}
\end{corollary}
\begin{proof}
$(ii)\Rightarrow (i)$ We show property $(iii)$ in Proposition \ref{lmn2}, i.e.\ that the map $\mathcal{I}(A)^ G \to \mathcal{I}(A\rtimes_r G ): I\mapsto I\rtimes_r G $ is surjective. Fix $J\in \mathcal{I}(A\rtimes_r G )$. Find a family of primitive ideals $(J_i)$ in $\mathcal{I}(A\rtimes_r G )$ together with irreducible representations $\pi_i\colon A\rtimes_r G \to B(H_i)$, such that
$$J=\cap_i J_i, \ \ \ \ker \pi_i=J_i$$
Set $I_i:=J_i\cap A$. From exactness of the action we get the canonical isomorphism, for every $i$, defined by
$$\iota_i\colon A\rtimes_r G /I_i\rtimes_r G  \to (A/I_i)\rtimes_r G :\ au_s + I_i\rtimes_r G  \mapsto (a+I_i)u_s,$$
and the well-defined (as $\pi_i(I_i\rtimes_r G )=0$) representation 
$$\pi^i\colon A\rtimes_r G /I_i\rtimes_r G  \to B(H_i): \ au_s + I_i\rtimes_r G  \mapsto \pi_i(au_s).$$
The map 
$$\pi^i\circ\iota_i^{-1}\colon (A/I_i)\rtimes_r G  \to B(H_i)$$
is {an irreducible} representation of $(A/I_i)\rtimes_r G $. With $J^{(i)}:=\ker \pi^i\circ\iota_i^{-1}$ we have that $J^{(i)}\cap A/I_i=0$ from the following argument: With $a+I_i\in J^{(i)}\cap A/I_i$
$$0=\pi^i\circ\iota_i^{-1}(a+I_i)=\pi^i(a+ I_i\rtimes_r G )=\pi_i(a).$$
Hence $a\in \ker\pi_i\cap A=I_i$ giving that $J^{(i)}\cap A/I_i=0$. Using $(ii)$ we obtain that $J^{(i)}=0$ and hence also $\ker\pi^i=0$. With
$$\rho_i\colon A\rtimes_r G \to A\rtimes_r G /I_i\rtimes_r G ,$$
fix $b\in J_i$. As $\pi^i(\rho_i(b))=\pi_i(b)=0$ it follows that $\rho_i(b)=0$. Hence 
$$J_i\subseteq \ker \rho_i=I_i\rtimes_r G \subseteq J_i, \ \ \ J=\cap_i J_i=\cap_i(I_i\rtimes_r G )=(\cap_i I_i)\rtimes_r G .$$
The last equality uses Proposition \ref{coo2}. The map $I\to I\rtimes_r G $ is surjective.
\end{proof}
\subsection{Exactness and essential freeness}
We will now present a way to ensure that $A$ separates the ideals in $A\rtimes_r G $ is by extending the result of Archbold and Spielberg in \cite{ArcSpi}. We let $\widehat{A}$ denote the spectrum of a C*-algebra $A$, i.e.\ the set of all equivalence classes of irreducible representations of $A$. The induced action of $ G $ on $\widehat{A}$ is given by $(t.x)(a):=x(t^{-1}.a), a\in A, [x]\in \widehat{A}$. For an action of $ G $ on a topological space $X$ we define the \emph{isotropy group} of $x\in X$ (also called the \emph{stabilizer subgroup}) as the set of all elements in $ G $ that fix $x$.
\begin{definition}[Boyle, Tomiyama,\mbox{ }\cite{BoyTom}]
\label{topfree}
{Let $G$ be a discrete group acting on a topological space $X$.} We say the action of $ G $ on $X$ is \emph{topologically free} provided that the points in $X$ with trivial isotropy are dense in $X$.
\end{definition}
\begin{definition}[Renault,\mbox{ }\cite{Ren:fixed}]
\label{essfree}
{Let $G$ be a discrete group acting on a topological space $X$.} We say the action of $ G $ on $X$ is \emph{essentially free} provided that for every closed invariant subset $Y\subseteq X$, the subset of points in $Y$ with trivial isotropy is dense in $Y$.
\end{definition}
\begin{remark}
\label{twotopfree}
Let $(A, G )$ be a C*-dynamical system with $ G $ discrete. Archbold and Spielberg define \emph{topological freeness} of $ G $ acting on $A$ slightly weaker than topological freeness of $ G $ acting on $\widehat{A}$. An action of $ G $ on $A$ is topologically free if $\cap_{t\in F}\{x\in \widehat{A}\colon t.x\neq x\}$ is dense in $\widehat{A}$ for any finite subset $F\subseteq G \setminus\{e\}$, cf.\ \cite{ArcSpi}. {If $\widehat{A}$ is Hausdorff (e.g. when $A$ is abelian) and $G$ 
is countable, then the two notions of topologically free agree because $\widehat{A}$ is always a Baire space.} If we use the natural inclusion $\widehat{A/I}\subseteq \widehat{A}$, then we see that \emph{the induced action of $ G $ on $A/I$ is topologically free for every $I\in \mathcal{I}(A)^G$} if an action of $ G $ on $\widehat{A}$ is essentially free. We will not make a weakening of the notion of essential
freeness, but we use sometimes the following concept:\\		
Let $\mathcal{P}$ denote a property for dynamical  systems  $(A,G)$. If this property holds (or is required) for all quotients $(A/I,G)$ with $I\in \mathcal{I}(A)^G$, then we say that $(A,G)$ is \emph{residually} $\mathcal{P}$. E.g.\ a conceptional name for the topological freeness for all quotients $A/I$ should be \emph{``residual'' topological	freeness} of the action on $A$.	
\end{remark}
Archbold and Spielberg considered when the reduced crossed product is simple. The key result was the following
\begin{theorem}[Archbold, Spielberg \mbox{}{\cite[Theorem 1]{ArcSpi}}]
\label{lmt1}
Let $(A, G )$ be a C*-dynamical system with $ G $ discrete and let $\pi$ be the surjection $A\rtimes G  \to A\rtimes_r G $. If the action of $ G $ on $A$ is topologically free then
$$\forall J\in \mathcal{I}(A\rtimes G ) : J\cap A=0 \Rightarrow \pi(J)=0.$$
\end{theorem}
Using the above result together with Theorem \ref{thetheorem} we get the following generalization.
\begin{theorem}
\label{tht1}
Let $(A, G )$ be a C*-dynamical system with $ G $ discrete. If the group $ G $ (or the action) is exact and the action of $ G $ on $\widehat{A}$ is essentially free then the algebra $A$ separates the ideals in $A\rtimes_r G $.
\end{theorem}
\begin{proof}
Fix $I\in \mathcal{I}(A)^ G $ and $J_1\in \mathcal{I}((A/I)\rtimes_r G )$ and assume $J_1\cap (A/I)=0$. Using Theorem \ref{thetheorem}.$(iii)$ it is enough to show that $J_1=0$. Let $\pi^{A/I}$ be the surjection $(A/I)\rtimes G  \to (A/I)\rtimes_r G $ and set $J:=(\pi^{A/I})^{-1}(J_1)$. Using that $\pi^{A/I}$ is just the identity on $A/I$ we get that
$$J\cap (A/I)=J_1\cap (A/I)=0$$
By Remark \ref{twotopfree} the action of $ G $ on $A/I$ is topologically free. Using Theorem \ref{lmt1} we obtain that $\pi^{A/I}(J)=0$. Hence $J_1=0$.
\end{proof}
\begin{remark}\label{rem:Gromov}
Gromov showed the existence of a finitely presented non-exact discrete group, cf.\ \cite{Gromov:exact}. Hence there exist a C*-dynamical system $(A, G )$ with a finitely presented discrete group  $G$  and a  non-exact action of $G$ on $A$. By Theorem \ref{thetheorem} $A$ does not separates the ideals in $A\rtimes_r G $.
	
If such $(A,G)$ can be found, such that the action is essentially free, then the Conjecture \ref{rencon} fails. But this is unknown if the action is essentially free, and the Conjecture \ref{rencon} remains an open problem. In fact, {these considerations} show that the Conjecture of Renault  is equivalent to the {question of whether essentially} free actions are exact.
\end{remark}
\begin{remark}
Let $(A, G )$ be a C*-dynamical system with $ G $ discrete. If the action is essentially free then the properties $(i)$--$(iii)$ in Proposition \ref{tho2} are all fulfilled. ({For a proof, use Theorem 1.15 for $A/I$, the automatic analogue of Lemma 1.9 for the full crossed product} and use that $(\rho\rtimes_r \id )\circ \pi^A=\pi^{A/I}\circ (\rho\rtimes \id )$ for the canonical maps $\rho\colon A\to A/I$, $\pi^B\colon B\rtimes G \to B\rtimes_r G $.)
\end{remark}
\begin{corollary}
\label{cot2}
Let $(A, G )$ be a C*-dynamical system with $ G $ discrete. Assume that the group $ G $ (or the action) is exact and that for every $x\in \widehat{A}$ the points in $\overline{ G .x}$ with trivial isotropy are dense in $\overline{ G .x}$. Then $A$ separates the ideals in $A\rtimes_r G $.
\end{corollary}
\begin{proof}
The action is exact and essentially free.
\end{proof}
\begin{corollary}[{cf.\mbox{ }\cite[Corollary 4.6]{Ren:fixed}}]
\label{cot3}
Let $(A, G )$ be a C*-dynamical with $ G $ discrete. Suppose the action is minimal (i.e.\ $A$ contains no non-trivial invariant ideals) and there exists an element in the spectrum of $A$ with trivial isotropy. Then $A\rtimes_r G $ is simple.
\end{corollary}
\begin{proof}
The action is exact and essentially free. 
\end{proof}
\begin{remark}
When considering the canonical action of $\Z$ on $\T=\R\cup\{\infty\}$ {(the right translation, fixing $\infty$)} it follows that topological freeness of $\Z$ acting on $\T$ is not enough to ensure that the action satisfy the residual intersection property. (For a proof consider 
{the short exact sequence
$$\xymatrix{0 \ar[r] & C_0(\R)\rtimes_r \Z\ar[r] & C(\T)\rtimes_r \Z \ar[r]  & \C\rtimes_r \Z\ar[r] & 0}.$$ If the action of $\Z$ on $C(\T)$  had the residual intersection property then the (trivial) action of $\Z$ on $\C$ would have the intersection property. But $\C\rtimes_r \Z \cong C^{*}(\Z ) \cong C(\T)$
has the property that all non-zero proper ideals (of which there are many) 
have zero intersection with the complex numbers $\C$.)}
\end{remark}
\begin{remark} In many cases, including when a countable group acts by an amenable action on a unital and abelian C*-algebra, essential freeness and the residual intersection property are equivalent. (For a proof examine Corollary \ref{lastresult}.) {For this class of examples, the essential freeness and exactness of the action are together necessary and sufficient to ensure that $A$ separates the ideals in $A\rtimes_r G $ (by Theorem 1.16 and Theorem 1.10).}
\end{remark}
\begin{remark}
Essential freeness and the residual intersection property are in general different conditions. (For a proof consider any simple crossed product where the action is not essentially free. For example the free group of two generators acting trivially on $\C$.)
\end{remark}
\subsection{The Rokhlin property}
One application of Theorem \ref{thetheorem} is an easy proof of the fact that the algebra separates the ideals in the reduced crossed product provided the action has the "Rokhlin property". 

In the following set $A_\infty:=l^\infty(A)/c_0(A)$, where $l^\infty(A)$ is the C*-algebra of all bounded functions from $\N$ into $A$ and where $c_0(A)$ is the ideal in $l^\infty(A)$ consisting of sequences $(a_n)_{n=1}^\infty$ for which $\|a_n\|\to 0$. We will suppress the canonical inclusion map $A\subseteq A_\infty$. The induced action of $ G $ on $A_\infty$ is defined entry-wise.
\begin{definition}[Izumi,\mbox{ }\cite{Izu:rohlin}]
\label{der1}
Let $(A, G )$ be a C*-dynamical system with $A$ unital and $ G $ finite. We say that the action (of $ G $ on $A$) satisfies the \emph{Rokhlin property} provided there exists a projection $p_e\in A'\cap A_{\infty}$ such that
\begin{itemize}
\item[$(i)$] $p_e\ \bot \ t.p_e, \ \ \ t\neq e$
\item[$(ii)$] $\sum_{t\in  G} t.p_e=1_{A_{\infty}}$
\end{itemize}
\end{definition}
\begin{theorem}
\label{thr1}
Let $(A, G )$ be a C*-dynamical system with $A$ unital and $ G $ finite. Suppose that the action satisfies the Rokhlin property. Then $A$ separates the ideals in $A\rtimes_r G $ .
\end{theorem}
\begin{proof}
Let $p_e\in A'\cap A_{\infty}$ be the projection ensuring the action satisfies the Rokhlin property and set $p_s=s.p_e$ for $s\in  G $. As $ G $ is finite the action is exact. Fix $x:=\sum_{t\in G }a_tu_t\in (A\rtimes_r G )^+$. Using Theorem \ref{thetheorem}.$(ii)$ it is enough to show that $E(x)\in \textrm{Ideal}[x]$. 

We can take the implementing unitaries for the action of $G$ on $A_\infty$ to be the same as those for $A$. We have the commuting triangle with three canonical inclusions
\begin{displaymath}
\xymatrix{A\rtimes_r G \ar[r] \ar[rd]& A_\infty\rtimes_r G \ar[d]\\
& (A\rtimes_r G)_\infty}
\end{displaymath}
giving the identity
$$\textrm{Ideal}_{A_\infty\rtimes_r G }[x]\cap A\rtimes_r G =\textrm{Ideal}_{A\rtimes_r G }[x].$$
{For the completeness of the proof let us recall how the equality is obtained. With $B:=A\rtimes_r G$ fix $b\in \textrm{Ideal}_{B_\infty}[x]\cap B$ and $\epsilon>0$. Denote the quotient map $l^\infty(B)\to B_\infty$ by $\pi_\infty$. Find $a=\sum_{j=1}^nt_jxs_j\in \textrm{Ideal}_{B_\infty}[x]$ such that $t_j=\pi_\infty(t_j^{(i)}), s_j=\pi_\infty(s_j^{(i)})$ and $\|a-b\|< \varepsilon$. From $\|\pi_\infty(\sum_{j=1}^nt_j^{(i)}xs_j^{(i)}-b)\|=\limsup_i\|\sum_{j=1}^nt_j^{(i)}xs_j^{(i)}-b\|<\varepsilon$ is follows that $b\in \textrm{Ideal}_{B}[x]$. Since $\textrm{Ideal}_{B}[x]\subseteq \textrm{Ideal}_{A_\infty\rtimes_r G}[x]\subseteq \textrm{Ideal}_{B_\infty}[x]$ and $\textrm{Ideal}_{B_\infty}[x]\cap B\subseteq \textrm{Ideal}_{B}[x]$ we are done.}

We only need to verify that $E(x)\in \textrm{Ideal}_{A_\infty\rtimes_r G }[x]$. This follows from the calculation below ($\delta_{e,t}$ is the Kronecker delta)
\begin{eqnarray*}
p_s(au_t)p_s&=&\delta_{e, t}p_s(au_t)=p_sE_{A}(au_t), \ \ \ s,t\in  G  \ a\in A,\\
\sum_{s\in G }p_sxp_s&=&\sum_{s\in G }p_sE_{A}(x)=E_A(x)\in \textrm{Ideal}_{A_\infty\rtimes_r G }[x].
\end{eqnarray*}
\end{proof}
\begin{remark}
Theorem \ref{thr1} only {applies} when the group $ G $ is finite. It is natural to consider if there is a way to extend this result. As we will see that is indeed the case.
\end{remark} 
\section{The residual Rokhlin* property}
The Rokhlin property is a quite restrictive property. We present here a weaker condition, called "residual Rokhlin* property", ensuring that $A$ separates the ideals in $A\rtimes_r G $, when considering a discrete exact group $ G $ acting on a C*-algebra $A$.

For a C*-dynamical system {$(A, G )$} the induced action of $ G $ on $A^*$ (resp. on $A^{**}$) is given by $(t.\varphi)(a):=\varphi(t^{-1}.a)$ for $a\in A, \varphi\in A^{*}$ (resp. for $a\in A^*, \varphi\in A^{**}$). We write $(A, G )\cong (B, G )$ when the isomorphism between $A$ and $B$ is equivariant, i.e.\ action preserving.
\begin{definition}	
\label{der2}
Let $(A, G )$ be a C*-dynamical system with $ G $ discrete. We say that the action (of $ G $ on $A$) has the \emph{Rokhlin* property} provided that there exist a projection $p_e\in A'\cap (A_\infty)^{**}$, such that
\begin{itemize}
\item[$(i)$] $p_e\ \bot \ t.p_e, \ \ \  t\neq e$
\item[$(ii)$] For every $a\in A$ with $a\neq 0$ there exist $t\in  G $ such that $a(t.p_e)\neq 0,$
\end{itemize}
If the induced action of $ G $ on $A/I$ has the Rokhlin* property for every invariant ideal $I$ in $A$, we say that {the action} (of $ G $ on $A$) has the \emph{residual Rokhlin* property}.
\end{definition}
\begin{remark}
Note that the residual Rokhlin* property is weaker than the Rokhlin property. This follows from the fact that the property $(ii)$ in Definition \ref{der2} is equivalent to the condition that $\|a\sum_{t\in G }t.p_e\|=\|a\|$ for every $a\in A$ {and the fact that the Rokhlin property automatically implies residual Rokhlin property. 
(Let us verify the second of the two mentioned properties. Suppose $p_e\in A'\cap A_\infty$ is the projection ensuring that the action of $G$ on $A$ has the Rokhlin property. Fix $I\in \mathcal{I}(A)^G$ and set $p_e^I=\varphi(p_e)$ using the canonical equivariant map $\varphi\colon A_\infty\to (A/I)_\infty$. We obtain that $p^I_e\ \bot \ t.p^I_e$ for $t\neq e$ and $\sum_{t\in  G} t.p^I_e=\sum_{t\in  G} \varphi(t.p_e)=\varphi(1_{A_{\infty}})=1_{(A/I)_{\infty}}$)}.
\end{remark}
\begin{example}
\label{example1}
Before giving {a proof of} {Theorem 2.5} we present a concrete C*-dynamical system $(A, G )$ with the residual Rokhlin* property. We consider the well-known example, where $A:=M_{n^\infty}\otimes {\mathcal K}$, the stabilized UHF-algebra of type $n^\infty$, where $ G :=\Z$, and where $\Z$ acts on $M_{n^\infty}\otimes {\mathcal K}$ via an automorphism $\bar\lambda$ that scales the trace {by a factor} $1/n$. In this case we have that $(M_{n^\infty}\otimes {\mathcal K})\rtimes_r\Z$ is isomorphic to ${\mathcal O}_n\otimes {\mathcal K}$, cf.\ \cite{Cuntz:On}.

More precisely let $(A,\mu_m\colon M_{n^\infty}\to A)$, cf.\ \cite[Definition 6.2.2]{RorLarLau:k-theory}, be the inductive limit of the sequence
\begin{displaymath}
\xymatrix{M_{n^\infty}\ar[r]^\lambda & M_{n^\infty} \ar[r]^\lambda & M_{n^\infty} \ar[r]^\lambda  & \cdots}
\end{displaymath}
with the connecting map $\lambda$ defined by $a\mapsto E_{11}\otimes a$, where $E_{11}$ is the projection $\textrm{diag}(1,0,\dots,0)$ in $M_n$. It is well known that $\lambda$ induces an automorphism $\bar\lambda$ on $A$ fulfilling that $\mu_m\circ\lambda=\bar\lambda\circ\mu_m$, cf.\ \cite[Prop. 6.2.4]{RorLarLau:k-theory}.
Using the identification $A\rtimes_r G \cong {\mathcal O}_n\otimes {\mathcal K}$ and the simplicity of $A$ we can give a new proof of the well known result, that the Cuntz algebra ${\mathcal O}_n$ is simple, by proving that the action of $ G $ on $A$ satisfies the Rokhlin* property. This follows using the projection
$$p_e:=(q_1,q_2,q_3\dots)\in A_\infty,$$
where $$q_k:=\mu_k(\underbrace{1\otimes\cdots\otimes 1}_{2k}\otimes (1-E_{11})\otimes \underbrace{E_{11}\otimes\cdots\otimes E_{11}}_{2k}\otimes 1\otimes \dots),$$ for $k\in\N$. Let us now show that $p_e\in A'\cap (A_\infty)^{**}$ and property $(ii)$. Fix $m, l\in\N$ and $a\in M_{n^m}$. Using the canonical inclusion $M_{n^m}\subseteq M_{n^\infty}$, together with the identities $M_{n^\infty}=\overline{\bigcup_{i=1}^\infty M_{n^i}}, A=\overline{\bigcup_{i=1}^\infty \mu_i(M_{n^\infty})}$, it is enough to show that $p_e\mu_l(a)=\mu_l(a)p_e$ and $\|\mu_l(a)\|=\|\mu_l(a)p_e\|$. But this follows from a simple calculation.
To see $(i)$ we simply use that the projections $(t.q_k)_{t\in\N}$ are pairwise orthogonal for a fixed $k\in\N$.
\end{example}
We now return to proving the fact that $A$ separates the ideals in $A\rtimes_r G $ provided a discrete exact group $ G $ acts on $A$ by an action with the residual Rokhlin* property. The key idea is the following Lemma
\begin{lemma}
\label{lmp2}
{Fix} a discrete group $ G $ and von Neumann algebras $N\subseteq M$. Suppose there exists an action of $ G $ on $N$ such that $(N, G )\cong (l_\infty( G ), G )$ (as usual $ G $ acts on $l_\infty( G )$ by left translation) and a group homomorphism $U\colon G \to \mathcal{U}(M)$ such that
$$U(t)fU(t)^*=t.f, \ \ \ t\in  G , f\in N.$$
Then $(B, G )$, given by
$$B:=N'\cap M, \ \ \ t.b:=U(t)bU(t)^*, \ \ \ t\in  G , b\in B,$$
is a C*-dynamical system. Further the representation
$$\id \times U\colon C_c( G ,B)\to M: bu_s\mapsto bU(s), \ \ \ s\in  G , b\in B$$
is isometric with respect to the reduced norm.
\end{lemma}
\begin{proof}
We show that $\id \times U$ is isometric. The first statement is straightforward. We want to use Fell's result \cite[Prop. 4.1.7]{BroOza}. Hence we make the identification $l_\infty( G )=N$ and $M\subseteq B(\HilH)$ for some Hilbert space $\HilH$.
Note that $(\id ,U,\HilH)$ is a normal covariant representation of $(l_\infty( G ), G )$ in the von Neumann algebra sense
(\footnote{A normal covariant representation of a discrete W*-dynamical system is a covariant representation in the usual C*-algebra sense with the additional assumption that the  representation of the von Neumann algebra is normal.
}).
It is a general fact, cf.\ \cite[Example 2.51]{Will:cross}, that any covariant representation $(\tilde{\pi},\tilde{U},\tilde\HilH)$ of $(l_\infty( G ), G )$, with $\tilde{\pi}\colon l_\infty( G )\to B(\tilde\HilH)$ unital and normal, is {unitarily} equivalent to a tensor multiple of the standard covariant representation 
$$l_\infty( G )\subseteq B(l_2( G )), \ \ \ t\mapsto \lambda_t\in B(l_2( G )),\ \ \ t\in  G ,$$ 
where the inclusion corresponds to multiplication operators $M_f, f\in l_\infty( G )$ and where $\lambda\colon G \to \mathcal{U}(l_2( G ))$ is the left regular representation. Thus, there exist a Hilbert space $H$ and a unitary $W\in B(l_2( G )\otimes H,\HilH)$, such that
$$W^*U(t)W=\lambda_t\otimes \id _H, \ \ \ W^*\id (f)W=M_f\otimes \id _H, \ \ \ t\in G, f\in l_\infty( G ).$$
With the canonical bijection $\varepsilon\colon B(l_2( G ))\bar\otimes B(H)\to B(l_2( G )\otimes H)$ and the unitary $V\in B(l_2( G )\otimes l_2( G )\otimes H)$, defined by $\delta_s\otimes \delta_t \otimes h\mapsto \delta_{t^{-1}s}\otimes \delta_t \otimes h$ where $(\delta_r)_{r\in G}$ is the basis for $l^2(G)$, consider the commuting diagram
\begin{displaymath}
\xymatrix{
(B, G ) \ar[rrr]^{(\id , U)} \ar[d]_{(\id , \id )} & & & M \ar[d]^{\id }\\
(l_\infty( G )'\cap B(\HilH), G ) \ar[rrr]^{(\id , U)} \ar[d]_{(W^*\cdot W, \id )} & & & B(\HilH) \ar[d]^{W^*\cdot W}\\
(\varepsilon(l_\infty( G )\bar\otimes B(H)), G ) \ar[rrr]^{(\id , \lambda\otimes 1)} \ar[d]_{(1\otimes\cdot, \id )} & & & B(l_2( G )\otimes H) \ar[d]^{\delta_s\otimes \delta_t \otimes h\mapsto \delta_s\otimes \cdot (\delta_t \otimes h)}\\
(1\otimes \varepsilon(l_\infty( G )\bar\otimes B(H)), G ) \ar[rrr]^{(1\otimes \id , 1\otimes \lambda\otimes 1)} \ar[d]_{(\id , \id )} & & & B(l_2( G )\otimes l_2( G )\otimes H) \ar[d]^{V^*\cdot V}\\
(1\otimes \varepsilon(l_\infty( G )\bar\otimes B(H)), G ) \ar[rrr]^{(1\otimes \id , \lambda\otimes \lambda\otimes 1)} & & & B(l_2( G )\otimes l_2( G )\otimes H)
}
\end{displaymath}
(For a C*-algebra $A\subseteq B(\mathscr K)$ the representation $1\otimes \id\colon 1\otimes A\to B(l^2(G)\otimes \mathscr K)$ sends $1\otimes a$ into the map defined by $\delta_s\otimes k\mapsto \delta_s\otimes a k$.) For completeness we show $$W^*(l_\infty( G )'\cap B(\HilH))W\subseteq \varepsilon(l_\infty( G )\bar\otimes B(H)).$$
{We have that $(l_\infty( G )\bar\otimes \C 1_{B(H)})'=l_\infty( G )'\bar\otimes B(H)=l_\infty( G )\bar\otimes B(H)$ where the first equality is a special case of the commutation theorem. Since both $\varepsilon$ and $W^*\cdot W$ map commutants to commutants it follows that
\begin{eqnarray*}
W^*(l_\infty( G )'\cap B(\HilH))W&=&\varepsilon\big(l_\infty( G )\bar\otimes \C 1_{B(H)}\big)'\\
&=&\varepsilon\big((l_\infty( G )\bar\otimes \C 1_{B(H)})'\big)=\varepsilon(l_\infty( G )\bar\otimes B(H)).
\end{eqnarray*}
}
Using Fell’s Absorption Principle II \cite[Prop. 4.1.7]{BroOza} the map $(1\otimes \id )\times (\lambda\otimes \lambda\otimes 1)$ is {unitarily equivalent to the regular representation} and hence is isometric with respect to the reduced norm $\|\cdot\|_\lambda$. 
From the commuting diagram it follows that $\id \times U$ is isometric with respect to $\|\cdot\|_\lambda$.\end{proof}
\begin{theorem}
\label{maintheorem} 
{ Let $(A, G )$ be a C*-dynamical system with $ G $ discrete.} Suppose the action is exact and {satisfies} the residual Rokhlin* property. Then $A$ separates the ideals in $A\rtimes_r G $.
\end{theorem}
\begin{proof}
By Theorem \ref{thetheorem} it is enough to show that the action of $ G $ on $A$ has the intersection property provided there exists a projection $p_e\in A'\cap (A_\infty)^{**}$ such that
\begin{itemize}
\item[$(i)$] $p_e\ \bot \ t.p_e, \ \ \  t\neq e$
\item[$(ii)$] For every $a\neq 0$ in $A$ there exist $t\in  G $ such that $a(t.p_e)\neq 0.$
\end{itemize}
From Remark \ref{lmp0} the intersection property is automatically fulfilled if every representation of $A\rtimes G $, which is faithful on $A$, weakly contains the regular representation of $(A, G )$. Fix a covariant representation $(\pi,u,H)$ of $(A, G )$ with $\pi$ faithful. As every representation of $A\rtimes  G $ comes from a covariant representation of $(A, G )$, cf.\ \cite{Dav:C*-ex}, it is sufficient to show that $\pi\times u$  weakly contains the regular representation, i.e.\ 
$$\|a\|_\lambda\leq \|\pi\times u(a)\|, \ \ \ a\in C_c( G ,A).$$
Using the commuting diagram
\begin{displaymath}
\xymatrix{&A\ \ar@{^{(}->}[r] \ar@{^{(}->}[d]^\pi & A_\infty\ \ar@{^{(}->}[r] \ar@{^{(}->}[d]^{\pi_\infty} & (A_\infty)^{**} \ar@{^{(}->}[d]^{ \pi_\infty^{**}} &\\
G  \ar[r]^{u}& B(H)\ \ar@{^{(}->}[r] & B(H)_{\infty}\ \ar@{^{(}->}[r] & (B(H)_\infty)^{**}\ \ar@{^{(}->}[r]& B(\HilH),}
\end{displaymath} 
{ where $(B(H)_\infty)^{**}$ is faithfully represented on a Hilbert space $\HilH$ and letting $u_\infty^{**}$ be the composition of the maps in the lower part of the diagram,} we get a covariant representation $(\pi_\infty^{**},u_\infty^{**},\HilH)$ of $((A_\infty)^{**}, G )$, with $\pi_\infty^{**}$ faithful, such that
\begin{eqnarray}
\label{eqnnorm}
\|\pi_\infty^{**}\times u_\infty^{**}(a)\|&=& \|\pi \times u(a)\|, \ \ \ a\in C_c( G ,A).
\end{eqnarray}
{Let $R$ be the weak$^∗$ closed algebra generated by $\{p_t\colon t\in G\}$ in $(A_\infty)^{**}$ with the identity $p$, where $p_t:=t.p_e$, $t\in G $. Define von Neumann algebras $N\subseteq M$ and an action of $ G $ on $N$ by
$$N:=\pi_\infty^{**}(R),\ \  M:=\pi_\infty^{**}(p)(B(H)_\infty)^{**}\pi_\infty^{**}(p),$$
$$s.\pi_\infty^{**}(a):=\pi_\infty^{**}(s.a),\ \ \ a\in R, \ \ \ s\in G.$$}
Note that $(N, G )\cong(l_\infty( G ), G )$, where $ G $ acts on $l_\infty( G )$ by left translation. With  
$$U(t):=\pi_\infty^{**}(p)u_\infty^{**}(t)\pi_\infty^{**}(p), \ \ \ t\in  G,$$
we get a group homomorphism $U\colon G \to \mathcal{U}(M)$ fulfilling that 
$$U(t)fU(t)^*=t.f, \ \ \ t\in  G , f\in N.$$ 
Set $B:=N'\cap M$. By Lemma \ref{lmp2}, the representation
$$\id \times U\colon C_c( G ,B)\to M: bu_s\mapsto bU(s), \ \ \ s\in  G , b\in B,$$
is isometric with respect to the reduced norm on $C_c( G ,B)$, and hence extends isometrically to 
$$\id \rtimes_r U\colon B\rtimes_r G \to M.$$ 
Now consider the homomorphism $\epsilon\colon A\to M, a\mapsto\pi_\infty^{**}(ap)$. {For every $a\in A$ we have that $ap$ commutes with the projection $p_t$ (since $(ap)p_t=ap_t=p_ta=p_t(pa)=p_t(ap)$) for $t\in G$. This implies that $ap$ commutes with the elements in $R$, and so $\epsilon(a)$ commutes with elements in $N$, i.e.\ $\epsilon(a)\in N'\cap M$. We conclude} that $\epsilon$ has image contained in $B$. 

By property $(ii)$ the map $a\mapsto ap$ on $A$ is faithful. As $\pi_\infty^{**}$ is faithful and the composition of faithful homomorphisms is faithful we conclude that $\epsilon$ is faithful. Furthermore $\epsilon\colon A\to B$ is equivariant and hence the canonical homomorphism $\epsilon \rtimes_r \id \colon A\rtimes_r G \to B\rtimes_r G $ is isometric.

When composing the two norm preserving maps $\id \rtimes_r U$ and $\epsilon \rtimes_r \id $ we obtain the following:
\begin{displaymath}
\xymatrix{
A \rtimes_r G   \ar[r]  & B \rtimes_r G  \ar[r] & M &&&\\
}
\end{displaymath}
\begin{displaymath}
\xymatrix{
\ \ \ \ \ \ & a  & \mapsto & \pi_\infty^{**}\times u_\infty^{**}(a)\pi_\infty^{**}(p) & a\in C_c( G ,A)
}
\end{displaymath}
It now follows that 
$$\|a\|_\lambda=\|\pi_\infty^{**}\times u_\infty^{**}(a)\pi_\infty^{**}(p)\|\leq \|\pi \times u(a)\|, \ \ \ a\in C_c( G ,A),$$
where the last inequality comes from (\ref{eqnnorm}). This completes the proof.
\end{proof}
\begin{remark}
The ideas in Theorem \ref{maintheorem} can be used to make more general results concerning C*-dynamical systems with locally compact groups.
\end{remark}	
\begin{corollary}
	\label{cor27}
{ Let $(A, G )$ be a C*-dynamical system with $ G $ discrete.} Suppose the action is exact and there exists a projection  $p_e\in A'\cap (A_\infty)^{**}$ such that
\begin{itemize}
\item[$(i)$] $p_e\ \bot \ t.p_e, \ \ \  t\neq e$
\item[$(ii)$] For every $a\in A$ and every closed invariant projection $q$ in the center of $A^{**}$ with $aq\neq 0$ there exist $t\in  G $ such that $aq(t.p_e)\neq 0$ \,
\end{itemize}
Then $A$ separates the ideals in $A\rtimes_r G $
\end{corollary}
\begin{proof}
It is a general fact that, for every C*-dynamical system $(B, G )$ with $ G $ discrete and $I\in \mathcal{I}(B)^ G $, we have the canonical equivariant inclusions $B\subseteq B_\infty$ and $B\subseteq B^{**}$. Further one has the decomposition $B^{**}=(B/I)^{**}+I^{**}$ with $(B/I)^{**}\cong B^{**}q_B^I$, where $q_B^I$ is the biggest central projection in $B^{**}$ orthogonal to $I$ ($q_B^I$ is the orthogonal complement to the supporting open central projection of $I$). Hence $q_B^I$ is a  closed invariant projection in the center of $B^{**}$.

Fix $I\in \mathcal{I}(A)^ G $ and set $B:=A/I$. Using the identification $B_\infty=A_\infty/I_\infty$ we have the commuting diagram
\begin{displaymath}
\xymatrix{A\ \ar@{^{(}->}[r] \ar[d] & A^{**}\ \ar@{^{(}->}[r] \ar[d]^{\cdot q_A^I} & (A_\infty)^{**} \ar[d]^{\cdot q_{A_\infty}^{I_\infty}} \\
B\ \ar@{^{(}->}[r] & B^{**}\ \ar@{^{(}->}[r] & (B_\infty)^{**}}
\end{displaymath}
Note that $q_A^I\leq q_{A_\infty}^{I_\infty}$ as elements in $(A_\infty)^{**}$. This follows from the fact that $q_A^I$ is central in $A^{**}$ (and hence also central in $(A_\infty)^{**})$ and orthogonal to $I$ (and hence also orthogonal to $I_\infty$).

Define $p_e^I:=p_eq_{A_\infty}^{I_\infty}\in B'\cap (B_\infty)^{**}$. {We now show that the action satisfies} the residual Rokhlin* property by verifying
\begin{itemize}
\item[$(i')$] $p_e^I\ \bot\ t.p_e^I, \ \ \  t\neq e$
\item[$(ii')$] For every $b\in B$ with $b\neq 0$ there exist $t\in  G $ such that $b(t.p_e^I)\neq 0.$
\end{itemize}
Property $(i')$ follows easily from $(i)$. To show $(ii')$ fix $b\in B$ with $b\neq 0$. By $(ii)$ and the left part of the diagram we have elements $a\in A$ and $t\in G $ such that $b=a+I$ and $aq_A^I(t.p_e)\neq 0$. From the right part of the diagram $b=aq_{A_\infty}^{I_\infty}$ in $(B_\infty)^{**}$. As $q_A^I\leq q_{A_\infty}^{I_\infty}$ we get that
\begin{eqnarray*}
b(t.p_e^I)&=&a(t.p_e)q_{A_\infty}^{I_\infty}=0 \textrm{ in }(B_\infty)^{**}\Rightarrow a(t.p_e)q_{A_\infty}^{I_\infty}=0 \textrm{ in }(A_\infty)^{**}\\
&\Rightarrow& aq_A^I(t.p_e)=(a(t.p_e)q_{A_\infty}^{I_\infty})q_A^I=0 \textrm{ in }(A_\infty)^{**},
\end{eqnarray*}
showing that $b(t.p_e^I)\neq 0$ and hence $(ii')$.
\end{proof}
\begin{remark}
	\label{better}
The condition $(ii)$ of Corollary \ref{cor27} can be replaced with a weaker condition without changing the conclusion. It is sufficient to consider only the closed invariant projections in the center of $A^{**}$ obtained as complement of support projections of invariant ideals in $A$.
\end{remark}
\subsection{The residual Rokhlin*-property and essential freeness}
We have shown in Theorem \ref{tht1} and Theorem \ref{maintheorem} that essential freeness and the residual Rokhlin*-property are sufficient to ensure that $A$ separates the ideals in $A\rtimes_r G $ for a discrete exact group $ G $ acting on $A$. 

We will now argue why the residual Rokhlin*-property is "better" than essential freeness. More precisely we will prove that in the case when $ G $ acts by an essentially free action on a C*-algebra $A$, then the residual Rokhlin* property is automatic. 
\begin{lemma}
\label{lmpt1b}
Let $(A, G )$ be a C*-dynamical system with $A$ abelian and $ G $ discrete. If $\varphi_1,\varphi_2$ are orthogonal states on $A$ and $\epsilon>0$ then there exist open sets $U_1,U_2\subseteq X:=\mathrm{Prim}(A)\cong\widehat{A}$ such that
$$U_1\cap U_2 =\emptyset, \ \ \ 1-\epsilon\leq \mu_i(U_i),$$
where $\mu_i$ is the unique regular Borel measure on $X$ corresponding to $\varphi_i$.
\end{lemma}
\begin{proof} By \cite[Definition 1.14.1]{Sak:C*-W*} we have that $\|\varphi_1-\varphi_2\|=2$. Find a selfadjoint element $h\in A$ such that $\|h\|\leq 1$ and $2-\epsilon\leq (\varphi_1-\varphi_2)(h)$. With 
$$h:=h_+-h_-,\ \ \  h_+,h_-\geq 0, \ \ \ h_+h_-=0$$ 
it follows that $1-\epsilon\leq \varphi_1(h_+)\leq \sup_n \varphi_1(h_+^{1/n})=\mu_1(\{x:h_+(x)>0\})$. {Similarly} we get that $1-\epsilon\leq \mu_2(\{x:h_-(x)>0\})$.
\end{proof}
\begin{lemma}
\label{lmpt2} Let $(A, G )$ be a C*-dynamical system with $A$ abelian, $ G $ discrete and let $\varphi$ be {a state} on $A$. {For every $t\in G $ such that $\varphi$ is orthogonal to $t.\varphi$, every $\epsilon > 0$ and} every open set $U\subseteq X:=\mathrm{Prim}(A)\cong\widehat{A}$ there exists an open set $U'\subseteq U$ such that
$$U'\cap t.U'=\emptyset,\ \ \ \mu(U')\geq\mu(U)-\epsilon,$$
where $\mu$ is the unique regular Borel measure on $X$ corresponding to $\varphi$.
\end{lemma}
\begin{proof} Assume having $t\in  G , \epsilon >0$ and $U\subseteq X$ open with the property that $\varphi$ is orthogonal to $t.\varphi$. Using Lemma \ref{lmpt1b} find open sets $U_1,U_2\subseteq X$ such that
$$U_1\cap U_2 =\emptyset, \ \ \ 1-\frac{\epsilon}{2}\leq \mu_i(U_i),$$
where $\mu_1:=\mu$ and $\mu_2:=t.\mu=\mu(t^{-1}.(\ \cdot \ ))$. Define
$$U'':=U_1\cap t^{-1}U_2, \ \ \ U':=U\cap U''.$$
Note that $U'\subseteq U$ is open and disjoint from $t.U'$. Using
\begin{eqnarray*}
\mu(U)&=&\mu(U\cap U'')+\mu(U\cap {U''}^c)\leq \mu(U')+\mu({U''}^c)\\
\mu({U''}^c)&=&\mu({(U_1\cap t^{-1}.U_2})^c)\leq \mu_1(U_1^c)+\mu_2(U_2^c)\leq \epsilon
\end{eqnarray*}
one get the desired property $\mu(U)\leq\mu(U')+\epsilon$ .
\end{proof}
\begin{theorem}
\label{tmpt1}
Let $(A, G )$ be a C*-dynamical system with $A$ abelian and $ G $ countable. Then the following are equivalent
\begin{itemize}
\item[$(i)$] The action of $ G $ on $\widehat{A}$ is topologically free.
\item[$(ii)$] There exists a projection $p_e\in A'\cap A^{**}${  ($=Z(A^{**})$)} such that
$$\|a\sum_{t\in G }(t.p_e)\|=\|a\|, \ \ \ p_e\ \bot \ s.p_e, \ \ \ s\in  G \setminus \{e\},\ a\in A.$$
\end{itemize}
In particular if the action is essentially free we obtain the residual Rokhlin* property.
\end{theorem}
\begin{proof}
$(i)\Rightarrow (ii)$:\, This direction works also in general for $(A,G)$ with discrete $G$ (and $A$ not necessary abelian and $G$ not necessary countable):
Let $F\subseteq \widehat{A}$ denote the set of points in $\widehat{A}$ with trivial isotropy. By axiom of choice we can find a subset $H\subseteq F$ such that $\{t.H\colon t\in  G \}$ is a partition of $F$. Recall that each von Neumann algebra canonically splits in the direct sum $N_{discrete}\oplus N_{continuous}$ by the maximal central projection $C_d$ such that $N_{discrete}:=NC_d$ is a discrete Type $I$ von Neumann algebra. In particular, for $N:=A^{**}$, the minimal projections of the center $Z(N_{discrete})$ of $N_{discrete}$ corresponds to the unitary equivalence classes of irreducible representation of $A$. In this way we get $G$-equivariant bijection between $\widehat{A}$ and the set of minimal projections in $Z(N_{discrete})$, i.e.\ {the center of} $A^{**}_{discrete}$ is naturally isomorphic to $l_\infty(\widehat{A})=c_0(\widehat{A})^{**}$, the bounded functions on $\widehat{A}$ equipped with the discrete topology. In particular, this gives   a one-to one correspondence between subsets $S\subseteq \widehat{A}$ of $\widehat{A}$ and projections $q_S$ in the center of $A^{**}$, such that $t.q_S=q_{t.S}$ for $t\in G$. The projection
$$p_e:=q_{H}$$
in the discrete part of $A^{**}$ satisfies the desired conditions giving $(ii)$, because the projection $\bigvee_{t\in G} t.q_H$ is the same as the projection $q_F$ corresponding to $F\subseteq \widehat{A}$, and $\| a q_F\| =\| a\|$ for all $a\in A$ (because the irreducible representations in $F$ are separating for $A$).

$(ii)\Rightarrow (i)$:\, 
Identify $A=C_0(X)$ for $X:=\mathrm{Prim}(A)\cong\widehat{A}$ and let $F\subseteq X$ be the subset of points in $X$ with trivial isotropy, i.e.\ elements that are {fixed only} by $e\in G $. 

Set $U:=X\backslash \bar{F}$ and $p_e^U:=p_eq$, where $q\in C_0(X)^{**}$ is the open invariant projection corresponding to the invariant ideal $C_0(U)$ in $C_0(X)$ such that 
$$C_0(U)^{**}\cong C_0(X)^{**}q, \ \ \ aq=a, \ \ \ a\in C_0(U).$$ 
Note that $U=\emptyset$ implies $(i)$. Assuming $U\neq \emptyset$ we show $U$ contains elements with trivial isotropy, giving a contradiction. First we show $p_e^U\neq 0$. Using that $A$ is weakly dense in $A^{**}$ find an increasing net $(a_i)$ of positive norm one elements in $C_0(U)$ weakly converging to $q\in A^{**}$. With $p:=\sum_{t\in G }t.p_e$ it follows from $(ii)$ that 
$$1= \|a_i\|=\|a_ip\|=\|a_iqp\|.$$
Since $p_e^U=0 \Rightarrow qp=0\Rightarrow a_iqp=0$, we get that $p_e^U\neq 0$. 

Find a normal state $\varphi^{**}$ on $C_0(U)^{**}$ with $\varphi^{**}(p_e^U)=1$ and let $\varphi$ be the restriction to $C_0(U)$. Now let $\mu$ be the regular Borel measure on $U$ corresponding to $\varphi$. Fix $\{t_1,t_2,\dots\}:= G \setminus \{e\}$ and $\epsilon_n=2^{-(n+1)}, n\in \N$. We now show, using induction, the existence of open sets $U_n\subseteq \dots\subseteq U_1\subseteq U$ fulfilling that 
$$U_n\cap t_j.U_n=\emptyset, \ j=1,\dots., n\ \ \ \mu(U_n)\geq 1-(\epsilon_1+\dots+\epsilon_n).$$
{For any $t\neq e$ the state $\varphi$ is orthogonal to $t.\varphi$ (since $\varphi^{**}(p_e^U)=1$ the projection $p_e^U$ is {larger} that the support projection of $\varphi$ and as $p_e^U\ \bot \ t.p_e^U$ the support projections of $\varphi$ and $t.\varphi$ are orthogonal).}

For $n=1$ we note that $\varphi$ is orthogonal to $t_1.\varphi$. Using Lemma \ref{lmpt2} on $t_1\in G $, $\epsilon_1>0$ and $U\subseteq X$ there is an open set $U_1\subseteq U$ such that
$$U_1\cap t_1.U_1=\emptyset,\ \ \ \mu(U_1)\geq\mu(U)-\epsilon_1=1-\epsilon_1.$$
For $n\geq 1$ we assume having the desired sets $U_1,\dots, U_n$. Using  Lemma \ref{lmpt2} on $t_{n+1}\in G  $ and $\epsilon_{n+1}>0$ and $U_n\subseteq X$ one gets an open set $U_{n+1}\subseteq U_n$ with
$$U_{n+1}\cap t_{n+1}.U_{n+1}=\emptyset,\ \ \ \mu(U_{n+1})\geq\mu(U_n)-\epsilon_{n+1}.$$
This implies the desired properties for $U_{n+1}$. As $\mu(\cap U_n)=\lim_n\mu(U_n)>0$ the set $\cap U_n$ is a non-empty subset in $U$ consisting of elements with trivial isotropy.
\end{proof}
{\begin{remark}
The property $(ii)$ in Theorem \ref{tmpt1} is different from either the Rokhlin or Rokhlin* statements. However it is clear that property $(ii)$ in Theorem \ref{tmpt1} implies the Rokhlin* property.
\end{remark}}
\begin{remark} 
\label{rem.ess.rokhlin}
We have the  following more general result:\\ 
{\it If $(A, G )$ is a C*-dynamical system  with $ G $ discrete and if the action of $ G $ on $\widehat{A}$ is essentially free, then $(A,G)$ satisfy the residual Rokhlin* property.}
See the proof of the direction $(i)\Rightarrow(ii)$ in the proof of Theorem \ref{tmpt1}.
\end{remark}
As a corollary we get a {particularly} nice reformulation of when $A$ separates the ideals in $A\rtimes_r G $ in the case where the action is amenable. The formulation of amenability stated below is taken from the book by Brown and Ozawa, cf.\  \cite{BroOza}.
\begin{definition}
{ Let $(C(X), G )$ be a C*-dynamical system with $X$ compact and $ G $ discrete.} The action is called \emph{amenable} if there exist a net of continuous maps $m_i\colon X\to \textnormal{Prob}( G )$ such that for each $s\in G $, $$\lim_{i\to\infty}\Big(\sup_{x\in X}\|s.m_i^x-m_i^{s.x}\|_1\Big)=0,$$ where $s.m_i^x(t)=m_i^x(s^{-1}t)$.
\end{definition}
\begin{corollary}
\label{lastresult}
{Let $(A, G )$ be a C*-dynamical system with $A$ unital, abelian and $ G $ countable, discrete}. Suppose the action is amenable. Then $A$ separates the ideals in $A\rtimes_r G $ if and only if the action satisfy the residual Rokhlin* property.
\end{corollary}
\begin{proof}
Assume $A$ separates the ideals in $A\rtimes_r G $. By Theorem \ref{thetheorem} the action of $ G $ of $A$ {satisfies} the residual intersection property. This can be reformulated in the following way. For every $I\in \mathcal{I}(A)^ G $ we have that
$$\forall J\in \mathcal{I}(A/I\rtimes G ) : \pi^{-1}(\pi(J))=J, J\cap A=0 \Rightarrow \pi(J)=0,$$
where $\pi$ is the canonical surjection $A/I\rtimes G \to A/I\rtimes_r G $. By \cite[Exercise 4.4.3]{BroOza} the induced action of $ G $ on $A/I$, $I\in \mathcal{I}(A)^ G $ is amenable implying that the map $\pi\colon A/I\rtimes G \to A/I\rtimes_r G $ is an isomorphism, cf.\  \cite[Theorem 4.2.6]{BroOza}. Using \cite[Theorem 2]{ArcSpi} the action of $ G $ on $A/I$ is topologically free for every $I\in \mathcal{I}(A)^ G $. {Since $G$ is countable, we obtain from Remark \ref{twotopfree} that the action of $G$ on $\widehat{A}$ is essentially free.} Using Theorem \ref{tmpt1} we get the residual Rokhlin* property.

Conversely assume that the action satisfy the residual Rokhlin* property. An amenable action is automatically exact. Hence $A$ separates the ideals in $A\rtimes_r  G $ by Theorem \ref{maintheorem}.
\end{proof}
\subsection{Additional remarks}
The following considerations show the connection between Rokhlin* property and
proper outerness.
\begin{definition}
In \cite{Ell} an automorphism $\alpha$ of $A$ is called \emph{properly outer} if for any non-zero \emph{$\alpha$-invariant} (i.e.\ $\langle \alpha \rangle$-invariant, $\langle \alpha \rangle$ being the group generated by $\alpha$) ideal $I$ in $A$ and any inner automorphism $\beta$ of $I$, $\|\alpha|_I-\beta\|=2$. 

{Let $(A, G )$ be a C*-dynamical system with $G$ discrete.} The action, say $\alpha$, of $G$ on $A$ is called \emph{properly outer} it $\alpha_t$ is properly outer for every $t\neq e$, cf. \cite{OlePed3}.
\end{definition}
Let $\mathcal{L}(A,A)$ be the set of linear and bounded maps from $A$ into $A$. We have a natural \emph{isometric and linear} map from $\mathcal{L}(A,A)$ to $\mathcal{L}(A_\infty, A_\infty)$ given by
$$\big(a\mapsto Ta\big) \to \big( (a_1,a_2,\cdots)+c_0(A) \mapsto (Ta_1,Ta_2,\cdots)+c_0(A)\big)$$
In particular, we get a natural unital monomorphism $\psi\colon \M(A) \to \M(A_\infty)$. The $\sigma((A_\infty)^{**}, (A_\infty)^*)$-closure of $A$ in $(A_\infty)^{**}$ is naturally W*-isomorphic to $A^{**}$ via a natural (not necessarily unital) monomorphism $\kappa\colon A^{**}\to (A_\infty)^{**}$.

For a unitary $u\in \M(A)$ we let $Ad(u)_\infty\colon A_\infty\to A_\infty$ denote the natural automorphism
$(a_1,a_2,\cdots)+c_0(A) \mapsto (ua_1u^*,ua_2u^*,\cdots)+c_0(A)$. Note that $Ad(u)_\infty$ is given by the inner automorphism $Ad(U)\colon A_\infty\to A_\infty, b\mapsto UbU^*$ for the unitary $U:=\psi(u)$. We identify $\M(A_\infty)$ naturally with its image in $(A_\infty)^{**}$ inducing $Ad(U)^{**}\colon(A_\infty)^{**}\to (A_\infty)^{**}, b\mapsto UbU^*$.
\begin{lemma}
\label{rem:A.2.inner.iso}
{ Fix} a unitary $u\in \M(A)$ and an automorphism $\alpha$ of $A$. With $q:=\kappa(1_{A^{**}})$, $U:=\psi(u)$ and $\beta:=\alpha_\infty$ we have the following identities
$$\beta^{**}(q)=q\ \ \textrm{ and } \ \ \mathrm{Ad}(U)^{**}(b)=b\ \textrm{for all} \ b\in (A'\cap (A_\infty)^{**})q.$$
\end{lemma}
\begin{proof}
The proof is a just an application of the two diagrams below.
\begin{displaymath}
\xymatrix{A^{**}\ \ar[r]^{\alpha^{**}} \ar[d]_{\kappa} & A^{**} \ar[d]_{\kappa} & \M(A) \ar[d]_{\psi} \ar[r]& A^{**} \ar[d]_{\kappa}\\
(A_\infty)^{**}\ \ar[r]^{\beta^{**}} & (A_\infty)^{**} & \M(A_\infty) \ar[r]^{\cdot q} & (A_\infty)^{**}
}
\end{displaymath}	
For an element $b\in A'\cap (A_\infty)^{**}=\kappa(A^{**})'\cap (A_\infty)^{**}$ we obtain that $$\mathrm{Ad}(U)^{**}(bq)=(Uq)b(Uq)^{*}=\kappa(u)b\kappa(u^*)=bq.$$
\end{proof}
\begin{lemma}
\label{rem:A.3.inner.iso}
Let $(A,G)$ be a C*-dynamical system with $G$ discrete. Assume the {action} of $G$ on $A$ has the Rokhlin* property. Then for every subgroup $H$ in $G$ and every $H$-invariant ideal $I$ in $A$ the restricted action of $H$ on $I$ has the Rokhlin* property.
\end{lemma}
\begin{proof} Using the equivalence relation $s \sim_H t \Leftrightarrow \exists h\in H \colon t=hs$ on $G$ let $F$ be a subset of $G$ with one element {from} each equivalent class. Further let $p$ be the supporting open central projection of $I$. Using the strong convergent sum $q_e:=\kappa(p)\sum_{t\in  F}p_t$ in $I'\cap (I_\infty)^{**}$ it follows that the action of $H$ on $I$ has the Rokhlin* property.
\end{proof}
\begin{theorem}	
\label{th:additional.remark}
Let $(A,G)$ be a C*-dynamical system with $G$ discrete. If the action has the Rokhlin* property then it is automatically properly outer.
\end{theorem}
\begin{proof}
Assume the action $\alpha$ of $G$ on $A$ is not properly outer. Find $t\neq e$, an $\langle \alpha_t \rangle$-invariant ideal $I$ in $A$ and $u\in \M(I)$ such that
$$\|\alpha_t|_I-Ad(u)\|<2.$$
By Lemma \ref{rem:A.3.inner.iso} it is sufficient to show that the action of $\langle \alpha_t \rangle$ on $I$ does not have the Rokhlin* property. Hence we can assume $G=\langle \alpha_t \rangle$ and $A=I$ and show that the action of $G$ on $A$ does not have the Rokhlin* property. With $\alpha:=\alpha_t$ we have an automorphism $\alpha$ of $A$ and a unitary $u\in \M(A)$  such that
$$\|\alpha - Ad(u)\|<2.$$
Define $q:=\kappa(1_{A^{**}})$, $U:=\psi(u)$ and $\beta:=\alpha_\infty$ as in Lemma \ref{rem:A.2.inner.iso}. The linear map $\beta^{**}-\mathrm{Ad}(U)^{**}$ is equal to the natural extension $(\alpha-\mathrm{Ad}(u))_\infty^{**}$ of $\alpha-\mathrm{Ad}(u)\in \mathcal{L} (A,A)$ to a bounded linear operator on $(A_\infty)^{**}$. In particular, we get $\| \beta^{**}-\mathrm{Ad}(U)^{**}\|<2$.

Assume there exist a projection $p_e\in A'\cap (A_\infty)^{**}$ such that $p_e\bot \beta^{**}(p_e)$. With $p:=p_eq$ it follows from Lemma \ref{rem:A.2.inner.iso} that
$$\mathrm{Ad}(U)^{**}(q)=q=\beta^{**}(q), \ \ \ \mathrm{Ad}(U)^{**}(p)=p.$$
Since $pq=qp=p$ we have that $(2p-q)^*(2p-q)=q$ and hence
$$2\| \beta^{**}(p) - p\| = \|\beta^{**} (2p-q)-(2p-q)\|<2.$$
Note that $p\bot \beta^{**}(p)$. This implies that $\|\beta^{**}(p) - p\|^2=\|\beta^{**}(p)+p\|\in \{0,1\}$. It follows that $p=0$. By definition of $q$ and Lemma \ref{rem:A.2.inner.iso} we obtain
$$a t.p_e=t.(p_et^{-1}.a)=t.(p_et^{-1}.(qa))=t.(p_eqt^{-1}.a)=0, \ \ \ a\in A, t\in G.$$
The action of $G$ on $A$ does not have the Rokhlin* property.
\end{proof} 	
\begin{corollary}
\label{cor.additional.remark}
{Let $(A, G )$ be a C*-dynamical system with $G$ discrete.} The properties 
\begin{itemize}
\item[(i)] The action of $G$ on $\widehat A$ is topologically free.
\item[(ii)] The action of $G$ on $A$ has the Rokhlin* property.
\item[(iii)] The action of $G$ on $A$ is properly outer.
\end{itemize}
fulfills the implications $(i)\Rightarrow (ii) \Rightarrow (iii)$. In addition if $A$ is abelian and $G$ countable then we obtain the implication $(iii) \Rightarrow (i)$ making all the conditions equivalent.
\end{corollary}
\begin{proof}
{ $(i)\Rightarrow (ii)$. We refer to the proof of Theorem \ref{tmpt1} ($(i)\Rightarrow (ii)$).

$(ii)\Rightarrow (iii)$. See Theorem \ref{th:additional.remark}.

$(iii)\Rightarrow (i)$. Suppose a countable discrete group $ G $ acts by $\alpha$ on an abelian C*-algebra $A$. Set $X:=\widehat A$ and $U_t:=\{x\in X\colon t.x\neq x\}$. If the action on $X$ is not topologically free there exist $t\neq e$ such that $U_t$ is not dense in $X$ (using that $ G $ is countable). Hence there exist an open non-empty subset $V$ in $U_t^c$. Note that $I:=C_0(V)$ is in a natural way a $\alpha_t$-invariant ideal in $A$ fulfilling that $t.x=x$ for all $x\in V$. In particular $\|\alpha_t|_I-id_I\|=0$. Hence the action is not properly outer.
}
\end{proof} 
\begin{remark}\label{rem:on.proper.outerness}
We know from the paper of Archbold and Spielberg \cite{ArcSpi} that topological freeness of $G$ on $\widehat A$ implies the action is properly outer, which is also contained in Corollary \ref{cor.additional.remark}.

If follows from the paper of Olesen and Pedersen \cite[Lemma 7.1]{OlePed3} (or of Kishimoto \cite{Kis:Auto}), that for a properly outer action on a \emph{separable} C*-algebra $A$ we obtain that for $b\in (A\rtimes G)^+$ and $\varepsilon >0$ there exist $x\in A^+$ with
$\| x\|\leq 1$, $\| xbx-xE(b)x\|<\varepsilon$ and $\| xE(b)x\| > \| E(b)\| -\varepsilon$.

An inspection of the proof of their Theorem 7.2 gives that the latter observation implies that any closed ideal $J$ of $A\rtimes G$ with $J\cap A=\{ 0\}$ must be contained in the kernel of $\pi\colon A\rtimes G \to A\rtimes _r G$ (in particular the intersection property holds).

We can apply this in a similar way to all quotients, and get the following result: \it{{Let $G$ be a discrete group acting by $\alpha$ on a separable C*-algebra $A$.} If $[\alpha_t]_I : A/I \to A/I$ is properly outer for every $ G $-invariant closed ideal $I\not=A$ and every $t\in  G  \setminus \{ e\}$, then we obtain the residual intersection property. Thus, if in addition the action of $G$ on $A$ is exact, then $A$ separates the ideals in $A\rtimes_r G$, {by Theorem 1.10.}}
\end{remark}

As a student of the phd-school OP-ALG-TOP-GEO the author is partially supported by the Danish Research Training Council. The research was partially completed while the author was at The Fields Institute for Research in Mathematical Sciences, Toronto.

\bibliographystyle{amsplain}
\providecommand{\bysame}{\leavevmode\hbox to3em{\hrulefill}\thinspace}
\providecommand{\MR}{\relax\ifhmode\unskip\space\fi MR }
\providecommand{\MRhref}[2]{%
\href{http://www.ams.org/mathscinet-getitem?mr=#1}{#2}
}
\providecommand{\href}[2]{#2}

\end{document}